\documentclass[leqno]{bcp73}

\usepackage{amsfonts}
\usepackage{graphicx}
\usepackage{amscd}
\usepackage{amsmath}
\usepackage{amssymb}

\usepackage{latexsym}
\usepackage{amsmath}
\usepackage{amssymb}
\usepackage{xypic}
\newtheorem{theorem}{Theorem}[section]
\newtheorem{proposition}[theorem]{Proposition}

\newtheorem{lemma}[theorem]{Lemma}
\newtheorem{definition}[theorem]{Definition}

\newcommand{\eins}{\mbox{$1\!\!\!\:{\rm I}$}}
\newcommand{\BF}{\mbox{$\mathcal{B}(\mathcal{F})$}}

\newcommand{\BH}{\mbox{$\mathcal{B}(\mathcal{H})$}}
\newcommand{\BK}{\mbox{$\mathcal{B}(\mathcal{K})$}}
\def\emm#1{{\it#1\/}}
\begin{document}

\title{Decompositions of Beurling type\\
for $E_0$-semigroups}
\author{Rolf Gohm}
\address{Department of Mathematics\\
University of Reading\\
Whiteknights, PO Box 220\\
Berkshire, RG6 6AX,
UK\\
E-mail: r.gohm@reading.ac.uk}
\date{}
\keywords{Beurling, $ E_0$-semigroup, unitary cocycle,
adapted, exact}
\abbrevauthors{R. Gohm}
\abbrevtitle{Decompositions of $\scriptstyle E_0$-semigroups}

\mathclass{Primary 46L55, 47D03; Secondary 81S25.}
\maketitlebcp

\abstract{We define tensor product decompositions of $E_0$-semigroups
with a structure analogous to a classical theorem
of Beurling. Such decompositions can be characterized by
adaptedness and exactness of unitary cocycles. For CCR-flows
we show that such cocycles are convergent.
}

\section*{Introduction.}
A well-known theorem of Beurling characterizes invariant subspaces
of the right shift on $\ell^2(\mathbb{N})$ by inner functions in the unit
disc. In this case the restriction of the right shift to a
nontrivial invariant subspace is automatically conjugate
(unitarily equivalent) to the original shift. This interesting
self-similar structure is fundamental in many respects. It is
the prototype of a very fruitful interaction between operator
theory and function theory, see for example \cite{Ni,FF}.

In this paper we want to study a somewhat analogous self-similar
structure for operators on a different level. While the original
setting concerns isometries and decompositions of the Hilbert
space into direct sums, we want to study $E_0$-semigroups, i.e.,
pointwise weak$^*$-continuous semigroups of unital
$*$-endomorphisms of $\BH$ for some complex separable Hilbert
space $\mathcal{H}$ (cf.\! \cite{Ar}),
and decompositions of the Hilbert space
into tensor products. To make the analogy visible, we present
in Section 1 the Nagy-Foia\c s functional model for $*$-stable
contractions and their characteristic functions in a suitable
way and in particular we emphasize a limit formula for the
characteristic function which is not made explicit in the
standard presentations. This analogy motivates the definition
of decompositions of Beurling type for $E_0$-semigroups in
Section 2. It is then shown that there is a reformulation in
terms of unitary cocycles for amplifications of the
$E_0$-semigroup. Relevant properties of the cocycles are
adaptedness and exactness.

In fact, if the $E_0$-semigroup is a CCR-flow (cf.\! \cite{Ar}),
then this leads to a setting which has been extensively studied by
quantum probabilists (cf.\! \cite{Pa}). We have a tensor product
of an initial Hilbert space with a symmetric Fock space over
$L^2[0,\infty)$, maybe with a multiplicity space $k$. Unitary
cocycles of amplifications of the quantized shift can be used to
construct quantum stochastic Markov processes in the form of
socalled Evans-Hudson flows. In many cases the cocycles can be defined
by quantum stochastic differential equations based on an integration
theory initiated by Hudson and Parthasarathy. Pioneering work in
establishing connections with $E_0$-semigroups has been done by
Bhat in \cite{Bh}. See also \cite{Li} for an up-to-date introduction
to the classification of various kinds of cocycles in this setting.

In this paper we do not adopt this differential point of view.
In fact, it is an interesting task to relate our arguments to it
and we shall do that elsewhere. The question posed by
decompositions of Beurling type concerns conjugacy between the
flow on the Fock space and the flow obtained by a cocycle
perturbation of its amplification. It seems that this question
has not been studied systematically. Here our analogy to the
Nagy-Foia\c s functional models turns out to be useful. In Section 3
we introduce the notion of a convergent unitary cocycle as a
sufficient criterion for exactness. The main result of Section 4
is that for CCR-flows we have a converse, so that
decompositions of Beurling type are characterized by the
occurrence of convergent cocycles. The analogy to Section~1
becomes very explicit at this point and we have a kind of
dictionary between these structures. In particular, the role
of the $*$-stable contraction is now played by a semigroup
of unital completely positive maps with an absorbing vector
state.

Similar results for single completely positive maps and
endomorphisms, i.e., with a discrete time parameter, have
been obtained by the author in \cite{Go}, Chapter 2, where
also some additional motivation for studying this kind of problems
can be found. As may be expected, we can often use the insight
from the discrete case but at some points we need to be more
careful to choose the correct cocycles. For CCR-flows much is
known about that and therefore our results are more complete in
this case. The success story of Beurling's theorem indicates,
in our view, that the study of decompositions of Beurling type
may be helpful in the process of deciphering the structure of
more general $E_0$-semigroups.


\section{Characteristic functions as limits.}
Let $\mathcal{H}$ be a complex separable Hilbert space and
$T: \mathcal{H} \rightarrow \mathcal{H}$ a contraction. We
have defect operators $D = \sqrt{\eins - T^* T}$ and
$D_* = \sqrt{\eins - T T^*}$ and defect spaces
$\mathcal{D} = \overline{D \mathcal{H}}$ and
$\mathcal{D}_* = \overline{D_* \mathcal{H}}$.
Then
\[
R = \bigg( \begin{array}{cc}
        T & D_* \\
        D & -T^* \\
        \end{array}
    \bigg): 
\mathcal{H} \oplus \mathcal{D}_* \rightarrow
\mathcal{H} \oplus \mathcal{D}
\]
is a unitary called the rotation matrix.

Using $H^2(\mathcal{D})$, the Hardy space of the unit disc with
values in $\mathcal{D}$, we can write the minimal isometric dilation
of $T$ as
$$
U: \mathcal{H} \oplus H^2(\mathcal{D})  \rightarrow 
\mathcal{H} \oplus H^2(\mathcal{D}), \ \quad
h \oplus f  \mapsto  Th \oplus (D h + z f),
$$
where $z$ is the complex variable so that multiplication by $z$
represents a one-sided right shift. The following factorization
of powers of $U$ is instructive. We denote by $H^2({\cal D}_{*,k})$ the space $H^2({\cal D})$ with ${\cal D}$ 
in all levels up to $k$ replaced by ${\cal D}_*$, in other words, if 
$\sum^\infty_{n=0} a_n z^n$ is in $H^2({\cal D}_{*,k})$ then $a_n \in {\cal 
D}_*$ for $n \leq k$ and $a_n \in {\cal D}$ for $n > k$. For $k\in\mathbb{N}_0$ define
\[
R_k: \mathcal{H} \oplus H^2(\mathcal{D}_{*,k}) \rightarrow
\mathcal{H} \oplus H^2(\mathcal{D}_{*,k-1}),
\]
\[
h \oplus \Big(\sum^\infty_{n=0} a_n z^n\Big)  \mapsto
(Th + D_* a_k) \oplus \Big( \sum^{k-1}_{n=0} a_n z^n + (D h - T^* a_k) z^k
+ \sum^\infty_{n=k+1} a_n z^n \Big),
\]
i.e., $R_k$ is a kind of leg-numbering notation indicating $R$
acting on $\mathcal{H}$ and on the $k$-th level on the right. Then
it is easy to check that for all $n\in\mathbb{N}$
\[
U^n (h \oplus f) = R_0 \ldots R_{n-1} (h \oplus z^n f)
\]
(where we identify $\mathcal{H} \oplus 0$ in $\mathcal{H} \oplus \mathcal{D}_*$
and $\mathcal{H} \oplus \mathcal{D}$)
\[
= T^n h \oplus \big( D T^{n-1} h + D T^{n-2} h z + \ldots
+ D h z^{n-1} + z^n f \big).
\]
A contraction is called $*$-stable if $T^{*n} h \to 0$ if $n\to\infty$
for all
$h \in \mathcal{H}$. Equivalent characterizations: (1) The minimal
isometric dilation $U$ of $T$ is (unitarily equivalent to) a
one-sided shift $S$. (2) $T^*$ is a restriction of $S^*$.
To show these (obviously sufficient) conditions for a
$*$-stable contraction $T$ one checks that
$$
C: \mathcal{H}  \rightarrow  H^2(\mathcal{D}_*), \ \quad
h  \mapsto  D_* (\eins - z T^*)^{-1} h,
$$
is an isometry which can be used to embed $\mathcal{H}$ into
$H^2(\mathcal{D}_*)$. Because $S^* C = C T^*$, where $S$ is
multiplication by $z$ in $H^2(\mathcal{D}_*)$, we have (2).
See \cite{FF}, Chapter IX, Theorem 6.4, where it is also shown
how to extend $C$ to a unitary $W$ between the two dilation
spaces intertwining $U$ and $S$, i.e., $W U = S W$.
Explicitly,
$$
W: \mathcal{H} \oplus H^2(\mathcal{D})  \rightarrow 
H^2(\mathcal{D}_*), \ \quad
h \oplus f  \mapsto  C h + \Theta_T f,
$$
where $\Theta_T(z) = -T + z D_* (\eins - z T^*)^{-1} D  |_{\mathcal{D}}$
gives an inner function $\Theta_T$ with values in the bounded
operators from $\mathcal{D}$ to $\mathcal{D}_*$. It is called the
characteristic function of $T$. The $S$-invariant subspace
$(C \mathcal{H})^\perp$ in $H^2(\mathcal{D}_*)$ can be written as
$\Theta_T H^2(\mathcal{D})$, which corresponds to an operator version
of Beurling's theorem. For $T$ we have a functional model on a
$S$-coinvariant subspace.
We want to call attention to the limit formula
\[
\hat{W} = \lim_{n\to\infty} R^*_{n-1} \ldots R^*_0.
\]
Note that $R^*_{n-1} \ldots R^*_0$ is nothing but the adjoint of
the product occurring in the formula for $U^n$. Here we have a limit
of unitaries 
and we assert that
this limit is an isometry $\hat{W}$ onto $0 \oplus H^2(\mathcal{D}_*)$
which corresponds to the unitary
$W: \mathcal{H} \oplus H^2(\mathcal{D}) \rightarrow
H^2(\mathcal{D}_*)$ introduced above.

It is not difficult to check this limit formula by direct
computation. In fact, by induction
\[
\displaylines{
R^*_{n-1} \ldots R^*_0 \Big(h \oplus \sum^\infty_{k=0} a_k z^k\Big)
= (T^{*n} h + T^{*(n-1)} D a_0 + \ldots + D a_{n-1})
\cr
\oplus\quad \Big[ (D_* h - T a_0) +
(D_* T^* h + D_* D a_0 - T a_1) z
\cr
+ \ldots + (D_* T^{*(n-1)} h + D_* T^{*(n-2)} D a_0 + \ldots
- T a_{n-1}) z^{n-1} + \sum^\infty_{k=n} a_k z^k \Big].
\cr}\]
The first summand of the orthogonal sum converges to zero 
as $n\to\infty$
because $T^{*n} \to 0$ strongly while the second summand converges to
$C h + \Theta_T \sum^\infty_{k=0} a_k z^k$.

We remark that there is a version of the theory of characteristic
functions for one-parameter semigroups of contractions which is
closely related to Lax-Phillips scattering theory
(cf.\! \cite{LP}, Chapter III).
We do not discuss it here because our 
only intention in this section
has been to provide an analogue for the following developments
and for this it is enough to consider the most elementary situation.

\section{Decompositions of Beurling type.}
Let $(\alpha_t)_{t\geq 0}$ be an $E_0$-semigroup, i.e.,
a pointwise weak$^*$-continuous semigroup of normal unital
$*$-endomorphisms of $\mathcal{B}(\tilde{\mathcal{H}})$, where
$\tilde{\mathcal{H}}$ is a complex separable Hilbert space.

\begin{definition}\rm
A decomposition $\tilde{\mathcal{H}} = \mathcal{H} \otimes \mathcal{K}$
is \emm{of Beurling type} with respect to $(\alpha_t)_{t\geq 0}$ if
$\eins \otimes \BK$ is invariant for $(\alpha_t)$ and if the
$E_0$-semigroup $(\gamma_t)_{t\geq 0}$ on $\BK$, given by
\[
\eins \otimes \gamma_t(y) = \alpha_t (\eins \otimes y)
\quad \mbox{for all } t\geq 0,\, y\in \BK,
\]
is conjugate to $(\alpha_t)$.
\end{definition}

As a motivation for this terminology the reader should compare
such decompositions with those in
Section 1 where the minimal isometric dilation of a
$*$-stable contraction is shown to yield an analogous structure
directly related to Beurling's theorem. Here we have a tensor
product instead of a direct sum and we consider mappings of a
higher level. Nevertheless, the analogy is useful as we shall see.

To construct and classify decompositions of Beurling type
it is necessary to make their structure more explicit. Conjugacy
of $(\alpha_t)$ and $(\gamma_t)$ means that there is a unitary
$w: \mathcal{H} \otimes \mathcal{K} \rightarrow \mathcal{K}$
such that
\[
\alpha_t(x) = w^* \gamma_t(w x w^*) w
\quad \mbox{for all } x \in \mathcal{B}
(\mathcal{H} \otimes \mathcal{K}).
\]
Let $\delta \in \mathcal{H}$ be a unit vector. Then we have an
embedding $\mathcal{K} \simeq \delta \otimes \mathcal{K}
\subset \mathcal{H} \otimes \mathcal{K}$. The orthogonal projection
from $\tilde{\mathcal{H}} = \mathcal{H} \otimes \mathcal{K}$ onto
$\mathcal{K}$ is denoted by $q$. By $\hat{w} \in \mathcal{B}
(\mathcal{H} \otimes \mathcal{K})$ we denote the isometry onto
$\mathcal{K}$ given by $w$, i.e., $\hat{w} \xi = \delta \otimes
w\xi$ for all $\xi \in \mathcal{H} \otimes \mathcal{K}$.
Then with $Id = Id_{\mathcal{H}}$ we get
\[
\alpha_t(x) = \hat{w}^* (Id \otimes \gamma_t)
(\hat{w} x \hat{w}^*) \hat{w}
= u_t (Id \otimes \gamma_t)(x) u^*_t,
\]
where $u_t = \hat{w}^*(Id \otimes \gamma_t)(\hat{w}) $.
The $E_0$-semigroup $(Id \otimes \gamma_t)_{t \geq 0}$ is called an
amplification of $(\gamma_t)_{t\geq 0}$, and it is easily checked
that $(u_t)_{t\geq 0}$ is a left unitary $(Id \otimes \gamma_t)$-cocycle,
i.e., a strongly continuous family $(u_t)_{t\geq 0}$ of unitaries on
$\mathcal{H} \otimes \mathcal{K}$ satisfying
\[
u_{t+s} = u_t (Id \otimes \gamma_t)(u_s)
\ \quad \mbox{for all }t,s \geq 0, \, u_0 = \eins.
\]
In fact, for every isometry $\hat{w} \in \mathcal{B}
(\mathcal{H} \otimes \mathcal{K})$ with
$\hat{w} \hat{w}^* = q$ the formula
$u_t = \hat{w}^*(Id \otimes \gamma_t)(\hat{w}) $
always defines such a cocycle. This expresses the fact that an
amplification $(Id \otimes \gamma_t)$ is always cocycle conjugate
to $(\gamma_t)$. See \cite{Ar}, 2.2.4 and 2.2.5, for more details.
In this paper, by `cocycle' we always mean a strongly continuous
left unitary cocycle, as above.

The second ingredient of a decomposition of Beurling type is the
fact that  $(\gamma_t)$ is obtained from  $(\alpha_t)$ by restriction.
We can write this property in terms of the cocycle. In fact, for
all $t\geq 0$, $ y\in \BK$ we have
\[
\alpha_t (\eins \otimes y) = u_t (Id \otimes \gamma_t)(\eins \otimes y) u^*_t
= u_t (\eins \otimes \gamma_t(y)) u^*_t
\]
and this always equals $\eins \otimes \gamma_t(y)$ if and only if
$u_t$ commutes with $\eins \otimes \gamma_t(\BK)$
(for all $t \geq 0$). The following definition helps to
summarize these observations.

\begin{definition}\rm
A $(Id \otimes \gamma_t)$-cocycle $(u_t)_{t\geq 0}$ is called
\emm{adapted} if $u_t$ commutes with $\eins \otimes \gamma_t(\BK)$
(for all $t \geq 0$). It is called \emm{exact} with respect to
$\mathcal{K}$ if there is an isometry $\hat{w} \in
\mathcal{B}(\mathcal{H} \otimes \mathcal{K})$ onto $\mathcal{K}$
such that
\[
u_t = \hat{w}^* \; (Id \otimes \gamma_t)(\hat{w})
\quad \mbox{for all } t\geq 0.
\]
\end{definition}

Our discussion above yields the following result.
\begin{theorem}
The decomposition $\mathcal{H} \otimes \mathcal{K}$ is
of Beurling type with respect to $(\alpha_t)$, such that
$\eins \otimes \gamma_t(y) = \alpha_t (\eins \otimes y)$
for all $t\geq 0$, $y\in \BK$, if and only if there
exists a $(Id \otimes \gamma_t)$-cocycle
$(u_t)$ which is adapted and exact with respect to
$\mathcal{K}$, such that
$\alpha_t(\cdot) = u_t (Id \otimes \gamma_t)(\cdot) u^*_t$.
\end{theorem}

We can reinterpret Theorem 2.3 in the following way.
To find all decompositions of Beur\-ling type for an $E_0$-semigroup,
start with (its conjugate version) $(\gamma_t)$
on $\BK$ and look for all adapted and exact
cocycles of its amplifications. Thus the emphasis of our
study shifts to cocycles.

Let us add some comments on adaptedness and exactness according
to Definition 2.2. Recall that a $(\beta_t)$-cocycle
$(u_t)_{t\geq 0}$ for an $E_0$-semigroup $(\beta_t)_{t\geq 0}$
is local in the sense of Powers \cite{Po}
or a gauge cocycle as in \cite{Ar}, 2.8, if $u_t$ commutes with
the range of $\beta_t$ (for all $t \geq 0$).
If in our setting
$\mathcal{H} = \mathbb{C}$ then we can identify $(Id \otimes \gamma_t)$
and $(\gamma_t)$, and in this case adapted cocycles are
nothing but gauge cocycles.
But for $dim \mathcal{H} \geq 2$ we have a
strict inclusion of $\eins \otimes \gamma_t(\BK)$ into
$(Id \otimes \gamma_t)(\mathcal{B}(\tilde{\mathcal{H}}))$ and thus an adapted
cocycle may fail to be a gauge cocycle. The term `adapted' comes from
stochastic processes, i.e., from the possibility to construct adapted
quantum stochastic processes by using such cocycles. See our introduction
and also \cite{Bh}, Chapter~9.

In a similar sense our notion of exactness with respect to
$\mathcal{K}$ is a slight generalization of the usual notion
of exactness which is obtained for $\mathcal{H} = \mathbb{C}$,
i.e., $q = \eins$ (cf.\! \cite{Ar}, 2.2.2).
Exactness does not depend on the choice of the
unit vector $\delta \in \mathcal{H}$. This can be seen by the
following reformulation of the definition which is clearly
equivalent to the original one and does not involve $\delta$.
If $w: \mathcal{H} \otimes \mathcal{K} \rightarrow \mathcal{K}$
is a unitary then define $(Id \otimes \gamma_t)(w):
\mathcal{H} \otimes \mathcal{K} \rightarrow \mathcal{K}$
to be the operator which maps $\xi \otimes \eta$ to
$\gamma_t(w_{\xi}) \eta$, where
$w_{\xi}: \mathcal{K} \rightarrow \mathcal{K},
\; \eta \mapsto w(\xi \otimes \eta)$. With this notation we
can say that $(u_t)$ is exact with respect to
$\mathcal{K}$ if there is a unitary
$w: \mathcal{H} \otimes \mathcal{K} \rightarrow \mathcal{K}$
such that $u_t = w^* \; (Id \otimes \gamma_t)(w)$ for all $t \geq 0$.

\section{Convergent cocycles.}
In this short section we discuss a simple analytic property
which is sufficient for a $(Id \otimes \gamma_t)$-cocycle
$(u_t)$ to be exact with respect to $\mathcal{K}$. We continue
to use the notation of the previous section.

\begin{definition}\rm
A $(\beta_t)$-cocycle $(u_t)$ for an $E_0$-semigroup $(\beta_t)$
is called
\emm{convergent} if $u^*_t$ converges for $t \to \infty$ in the
strong operator topology.
\end{definition}

\begin{proposition}
If a $(Id \otimes \gamma_t)$-cocycle $(u_t)$ is
convergent and $u^*_t \to \hat{w}$ with $\hat{w}\hat{w}^* = q$,
then $(u_t)$ is exact with respect to $\mathcal{K}$.
\end{proposition}
\noindent
\proof
From the cocycle equation
\[
u_t = u_{t+s}(Id \otimes \gamma_t)(u^*_s),
\]
which for $s \to \infty$ converges weakly to
$\hat{w}^* (Id \otimes \gamma_t)(\hat{w})$, we conclude that
$u_t  =\break  \hat{w}^*  (Id \otimes \gamma_t)(\hat{w})$
for all $t\geq 0$, i.e., $(u_t)$ is exact with respect to $\mathcal{K}$.
\endproof

Note that, for a convergent cocycle, $(u_t)$ always converges
in the weak operator topology. Recall that in the unitary
group strong and weak convergence coincide. But in our
applications typically the limit is not unitary and $(u_t)$
does not converge strongly. In this respect the following
observation is useful.

\begin{lemma}
Let $(x_t)_{t \geq 0}$ be unitaries converging to $x$
in the strong operator topology. Then a vector $\xi$
is in the range of $x$ if and only if $\big(x^*_t(\xi)\big)$
is convergent in norm.
\end{lemma}

The lemma can be proved by combining the following elementary facts:
$(x^*_t)$ converges to $x^*$ in the weak operator topology.
Weakly convergent unit vectors are convergent in norm if
and only if the limit is a unit vector.
$\xi$ is in the range of $x$ if and only if $\|x^*(\xi)\| = \|\xi\|$.

\section{CCR-flows.}
The most basic examples of (non-automorphic) $E_0$-semigroups are
the CCR-flows and they can be realized on symmetric Fock spaces.
See \cite{Pa} for details on Fock space and \cite{Ar} for a
presentation of CCR-flows from the point of view of $E_0$-semigroups.
We want to discuss decompositions of Beurling type in this case.

Let us introduce some notation. We denote by
$\mathcal{F} := \Gamma(L^2[0,\infty),k)$ the symmetric Fock space over
$L^2[0,\infty)$ with a multiplicity space $k$. For all $t \geq 0$
we have a canonical decomposition
$\mathcal{F} = \mathcal{F}_{t]} \otimes \mathcal{F}_{[t}$, where
$\mathcal{F}_{t]} := \Gamma(L^2[0,t),k)$ and
$\mathcal{F}_{[t} := \Gamma(L^2[t,\infty),k)$.
There is also a canonical way to identify $\mathcal{F}$ and
$\mathcal{F}_{[t}$ using the right shift. Always suppressing the
notation for the identification maps we can define
\[
\gamma_t: \mathcal{B}(\mathcal{F}) \rightarrow \mathcal{B}(\mathcal{F}),
\ \quad
x \mapsto \eins_{t]} \otimes x,
\]
where $\eins_{t]}$ is the identity on $\mathcal{F}_{t]}$.
The $E_0$-semigroup $(\gamma_t)_{t \geq 0}$ is called the CCR-flow
of index $\dim  k$. In \cite{Po}, Powers gave the following abstract
characterization: An $E_0$-semigroup which is completely spatial
and in standard form is conjugate to a CCR-flow. See \cite{Ar},
2.6 and 2.7, for complete spatiality. `Standard form' means that
there exists an absorbing vector state. For a CCR-flow it is given
by the vacuum vector $\Omega \in \mathcal{F}$, and the absorbing
property means that
\[
\rho \circ \gamma_t \to \langle \Omega, \cdot\, \Omega \rangle
\quad (t \to \infty)
\]
for all normal states $\rho$. Note that an absorbing state is
always a unique invariant state.

Now let $\mathcal{H}$ be another Hilbert space and assume that
the $E_0$-semigroup $(\alpha_t)_{t \geq 0}$ on
$\mathcal{B}(\mathcal{H} \otimes \mathcal{F})$ provides a
decomposition of Beurling type such that
$\eins \otimes \gamma_t(y) = \alpha_t (\eins \otimes y)$
for all $t\geq 0$, $y\in \BF$. In particular,
$(\alpha_t)$ and $(\gamma_t)$ are conjugate and there is
also an absorbing vector state given by $\tilde{\Omega}
\in \mathcal{H} \otimes \mathcal{F}$ for $(\alpha_t)$.

\begin{lemma}
There is a unit vector $\delta \in \mathcal{H}$ such that
$\tilde{\Omega} = \delta \otimes \Omega$.
\end{lemma}
\noindent
\proof
For all $y\in \BF$ we obtain
\[
\langle \tilde{\Omega}, (\eins \otimes y) \tilde{\Omega} \rangle
= \langle \tilde{\Omega}, \alpha_t (\eins \otimes y) \tilde{\Omega} \rangle
= \langle \tilde{\Omega}, (\eins \otimes \gamma_t)(y) \tilde{\Omega} \rangle
\stackrel{t \to \infty}{\longrightarrow}
\langle \Omega, y\Omega \rangle.
\]
We conclude that $\langle \tilde{\Omega}, 
(\eins \otimes y) \tilde{\Omega} \rangle
= \langle \Omega, y\Omega \rangle$ for all $y\in \BF$. This means that
$\langle \tilde{\Omega}, \cdot\, \tilde{\Omega} \rangle$ is not an
entangled state. In other words, there is a unit vector
$\delta \in \mathcal{H}$ such that
$\tilde{\Omega} = \delta \otimes \Omega$.
\endproof

From Lemma 4.1 we see that here we have a canonical choice for the unit
vector $\delta \in \mathcal{H}$ which was chosen arbitrarily in Section 2.
This particular choice will be used in the following, together with the
canonical embeddings
\[
\mathcal{H} \simeq \mathcal{H} \otimes \Omega \subset
\mathcal{H} \otimes \mathcal{F} = \tilde{\mathcal{H}},
\]
\[
\mathcal{F} \simeq \delta \otimes \mathcal{F} \subset
\mathcal{H} \otimes \mathcal{F} = \tilde{\mathcal{H}}.
\]
We denote by $p = \eins_{\mathcal{H}} \otimes |\Omega\rangle \langle\Omega|$
the orthogonal projection with range $\mathcal{H}$ and by
$q = |\delta\rangle \langle\delta| \otimes  \eins_{\mathcal{F}}$
the orthogonal projection with range $\mathcal{F}$.

Now we can give some equivalent characterizations of decompositions
of Beurling type for CCR-flows. Note that by Theorem 2.3 the
presentation $\alpha_t(\cdot) = u_t (Id \otimes \gamma_t)(\cdot) u^*_t$
with an adapted cocycle $(u_t)$ is not a restriction of generality.
Let us call
the cocycle $(\hat{u}_t)$ a modification of $(u_t)$ if there is a
gauge cocycle $(\hat{v}_t)$ such that $\hat{u}_t = u_t  \hat{v}_t$
for all $t \geq 0$. Then also
$\alpha_t(\cdot) = \hat{u}_t (Id \otimes \gamma_t)(\cdot) \hat{u}^*_t$.
Further we need the compression of $\alpha_t$ to $\mathcal{H}$, i.e.,
\[
Z_t: \BH \rightarrow \BH,\quad x \mapsto p\alpha_t(pxp)p
= p\alpha_t(xp)p ,
\]
which is a unital completely positive map (for each $t \geq 0$).
\begin{theorem}
Let $(\gamma_t)$ be the CCR-flow on $\mathcal{F}$ and $(\alpha_t)$
another $E_0$-semigroup on $\mathcal{H} \otimes \mathcal{F}$
given by $\alpha_t(\cdot) = u_t (Id \otimes \gamma_t)(\cdot) u^*_t$
with an adapted $(Id \otimes \gamma_t)$-cocycle $(u_t)$.
The following assertions are equivalent:
\begin{itemize}
\item[\rm(a)]
$\tilde{\mathcal{H}} = \mathcal{H} \otimes \mathcal{F}$ is a
decomposition of Beurling type with respect to $(\alpha_t)$, such that
$\eins \otimes \gamma_t(y) = \alpha_t (\eins \otimes y)$
for all $t\geq 0,\, y\in \BF$.
\item[\rm(b)]
There exists a modification $(\hat{u}_t)$ of $(u_t)$ which is
exact with respect to $\mathcal{F}$, i.e.,
$\hat{u}_t = \hat{w}^* \; (Id \otimes \gamma_t)(\hat{w})$,
where $\hat{w}\,\hat{w}^* = q$.
\item[\rm(c)]
There exists a modification $(\hat{u}_t)$ of $(u_t)$ which is
convergent, such that $\hat{u}^*_t \to \hat{w}$,
where $\hat{w}\,\hat{w}^* = q$.
\item[\rm(d)]
$(Z_t)_{t\geq 0}$ is a unital CP-semigroup with an absorbing
vector state (given by $\delta \in \mathcal{H}$).
\item[\rm(e)]
$(Z_t)_{t\geq 0}$ is an ergodic unital CP-semigroup with an invariant
vector state (given by $\delta \in \mathcal{H}$).
`Ergodic' means that the only fixed points for all $Z_t$ are
scalar multiples of~$\eins_{\mathcal{H}}$.
\end{itemize}
\end{theorem}

We see that under the conditions of Theorem 4.2 we have a converse
of Proposition~3.2. Moreover, the analogy to the Hilbert
space theory in Section 1 becomes more explicit. We can think of
$(\alpha_t)$ as a homomorphic dilation of a semigroup $(Z_t)$
of unital CP-maps and this semigroup already determines whether
we have a decomposition of Beurling type. Moreover we have
[in part (c)] a limit formula which is parallel to the one
discussed in Section 1 and we may think of $\hat{w} \in
\mathcal{B}(\mathcal{H} \otimes \mathcal{F})$
or the corresponding $w: \mathcal{H} \otimes \mathcal{F}
\rightarrow \mathcal{F}$ as an analogue of the characteristic
function.

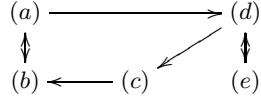
\proof
\[
\xymatrix@R=11pt{
(a) \ar[rr] \ar[d] & & (d) \ar[dl] \ar[d] \\
(b) \ar[u] & \ar[l] (c) & (e) \ar[u]
}
\]
The equivalence (a)$\Leftrightarrow$(b) is a reformulation of
Theorem 2.3 which we have repeated here for convenience.
(c)$\Rightarrow$(b) is Proposition 3.2.

Let us now show (a)$\Rightarrow$(d). The semigroup property of $(Z_t)$
follows from
\[
Z_t Z_s(x) = p \alpha_t(p\alpha_s(x)p)p
= p\alpha_t(p)\alpha_{t+s}(x)\alpha_t(p)p
\]
together with
\[
\alpha_t(p) = u_t (Id \otimes \gamma_t)(p)u^*_t
= u_t (Id \otimes \gamma_t)
(\eins_{\mathcal{H}} \otimes |\Omega\rangle \langle\Omega|)u^*_t
\]
\[
= u_t (\eins_{\mathcal{H}} \otimes
\gamma_t(|\Omega\rangle \langle\Omega|))u^*_t
= \eins_{\mathcal{H}} \otimes \gamma_t(|\Omega\rangle \langle\Omega|)
= \eins_{\mathcal{H}} \otimes \eins_{t]} \otimes |\Omega\rangle \langle\Omega|
\geq p,
\]
where we used the adaptedness of the cocycle $(u_t)$.
By Lemma 4.1 we know that $\delta \otimes \Omega$ yields an
absorbing vector state for $(\alpha_t)$. Thus for any
normal state $\rho$ on $\BH$
\[
\rho \circ Z_t(x) = (\rho \otimes \langle \Omega, \cdot\,\Omega \rangle)
\circ \alpha_t(xp) \stackrel{t \to \infty}{\longrightarrow}
\langle \delta \otimes \Omega, (xp) \delta \otimes \Omega \rangle
= \langle \delta,x\delta \rangle,
\]
i.e., $\langle \delta, \cdot\delta \rangle$ is an absorbing
vector state for $(Z_t)$.

(d)$\Leftrightarrow$(e) is well known, see for example \cite{Go}, A.5.2.
For convenience we sketch briefly the main arguments. For a nontrivial fixed
point we can find two normal states assigning different values to it
which clearly rules out the absorption property. Conversely, start
with the observation that for an invariant state
$\langle \delta, \cdot\,\delta \rangle$ we always have
$|\delta\rangle \langle \delta| \leq Z_t (|\delta\rangle \langle \delta|)$
and thus there is a strong limit of $(Z_t (|\delta\rangle \langle \delta|))$
for $t \to \infty$. Because this limit is a fixed point, ergodicity
implies that it coincides with $\eins_{\mathcal{H}}$. From this it follows
that the invariant state is absorbing.

It remains to prove (d)$\Rightarrow$(c). We start by showing that
there is a modification $(\hat{u}_t)$ of $(u_t)$ such that
\[
\hat{u}_t \; \delta \otimes \Omega = \delta \otimes \Omega
\quad {\mbox for\; all}\; t \geq 0.
\]
We have for all $t \geq 0$ and $x \in \BH$
\[
Z_t(x) = pu_t(Id \otimes \gamma_t)(xp) u^*_tp
= pu_t(x\otimes\eins)u^*_tp
\]
\big[using the definition of $(\gamma_t)$ and adaptedness of $(u_t)$\big]
and thus
\[
\langle u^*_t \delta \otimes \Omega,
(x\otimes\eins) u^*_t \delta \otimes \Omega \rangle
= \langle \delta \otimes \Omega, u_t(x\otimes\eins) u^*_t
\delta \otimes \Omega \rangle
= \langle \delta, Z_t(x) \delta \rangle
= \langle \delta, x \delta \rangle.
\]
We conclude that for all $t \geq 0$ we have
\[
u^*_t\; \delta \otimes \Omega = \delta \otimes \hat{\Omega}_t,
\]
where $\hat{\Omega}_t \in \mathcal{F}_{t]}$ is a unit vector.
The cocycle property of $(u_t)$ implies that
\[
\hat{\Omega}_t \otimes \hat{\Omega}_s = \hat{\Omega}_{t+s}
\]
(using the canonical identifications with respect to
$\mathcal{F}_{t]} \otimes \mathcal{F}_{s]} = \mathcal{F}_{t+s]}$).
This means that we can think of $(\hat{\Omega}_t)_{t \geq 0}$
as a unit. Explicitly, in Arveson's sense (cf.\! \cite{Ar}, 2.5)
the semigroup $(s_t)_{t \geq 0}$ given by the isometries
\[
s_t: \mathcal{F} \rightarrow \mathcal{F}, \quad
\eta \mapsto \hat{\Omega}_t \otimes \eta
\]
is a unit for $(\gamma_t)$.
It is known that the gauge cocycles act transitively on this
kind of units, see \cite{Ar},~3.8.5. In particular, we can
find a gauge cocycle $(v_t)_{t \geq 0}$ of $(\gamma_t)$
which maps the distinguished unit $(\Omega_t)_{t \geq 0}$ obtained by
factorization of the vacuum vector $\Omega$ to the
unit $(\hat{\Omega}_t)_{t \geq 0}$. Now we define
the modification $(\hat{u}_t)$ of $(u_t)$ by
\[
\hat{u}_t = u_t \; (\eins \otimes v_t)
\quad {\mbox for\; all}\; t \geq 0.
\]
Then $\hat{u}_t \; \delta \otimes \Omega = \delta \otimes \Omega$,
as claimed.

Now because $\langle \delta, \cdot\, \delta\rangle$ is absorbing for
$(Z_t)$, we can start with an arbitrary unit vector $\xi\in
\mathcal{H}$, and for all $x\in\BH$ we get
\[
\langle \hat{u}^*_t \xi \otimes \Omega,
(x\otimes\eins) \hat{u}^*_t \xi \otimes \Omega \rangle
=  \langle \xi, Z_t(x) \xi \rangle
\stackrel{t \to \infty}{\longrightarrow}
\langle \delta, x \delta \rangle.
\]
If we choose for $x$ the one-dimensional projection
$|\delta \rangle \langle \delta|$ then it follows that
\[
\| \big(|\delta \rangle \langle \delta| \otimes \eins\big) \;
\hat{u}^*_t \xi \otimes \Omega \|^2
\stackrel{t \to \infty}{\longrightarrow}
\langle \delta, \delta \rangle \langle \delta, \delta \rangle = 1
\]
and because $|\delta \rangle \langle \delta| \otimes \eins = q$
we have
\[
\| q \hat{u}^*_t \xi \otimes \Omega \|
\stackrel{t \to \infty}{\longrightarrow} 1.
\]
This means that we can choose vectors $\eta_t \in \mathcal{F}_{t]}$
such that
\[
\Delta_t := \| \hat{u}^*_t \xi \otimes \Omega - \delta \otimes \eta_t \|
\stackrel{t \to \infty}{\longrightarrow} 0.
\]
Then for $t,s \geq 0$
\[
\| (\hat{u}^*_{t+s} - \hat{u}^*_t) \xi \otimes \Omega \|
= \| \big[ (Id \otimes \gamma_t)(\hat{u}^*_s)\hat{u}^*_t - \hat{u}^*_t \big]
\xi \otimes \Omega \|
\]
\[
\leq 2\Delta_t + \| (Id \otimes \gamma_t)(\hat{u}^*_s)\delta \otimes \eta_t
- \delta \otimes \eta_t \| = 2\Delta_t.
\]
Indeed, $(Id \otimes \gamma_t)(\hat{u}^*_s)\delta \otimes \eta_t
= \delta \otimes \eta_t$ because $\eta_t \in \mathcal{F}_{t]}$
and $\hat{u}^*_s \; \delta \otimes \Omega = \delta \otimes \Omega$.
We conclude that $(\hat{u}^*_t\; \xi \otimes \Omega)_{t \geq 0}$
is a Cauchy net and thus convergent.

We have to prove that $(\hat{u}^*_t\; \hat{\xi})$ is convergent
for all $\hat{\xi} \in \mathcal{H} \otimes \mathcal{F}$. But this
can be reduced to the special case above as follows. First, it is
enough to consider $\hat{\xi} = \xi^\prime \otimes \eta^\prime_t$
with $\eta^\prime_t \in \mathcal{F}_{t]}$ for all $t$, because
linear combinations of such vectors are dense in
$\mathcal{H} \otimes \mathcal{F}$. Second, we have
\[
\hat{u}^*_{t+s}\; (\xi^\prime \otimes \eta^\prime_t)
= (Id \otimes \gamma_t)(\hat{u}^*_s)\;\hat{u}^*_t
\; (\xi^\prime \otimes \eta^\prime_t),
\]
and because $\hat{u}^*_t (\xi^\prime \otimes \eta^\prime_t)$
can be approximated by linear combinations of vectors of the form
$\xi \otimes \eta_t$ with $\eta_t \in \mathcal{F}_{t]}$
we see that convergence of
$\hat{u}^*_{t+s} (\xi^\prime \otimes \eta^\prime_t)$ for
$s \to \infty$ always takes place if $\hat{u}^*_s (\xi \otimes \Omega)$
converges for all $\xi \in \mathcal{H}$. But the latter has already
been shown.

We have established the strong convergence of $(\hat{u}^*_t)$, i.e.,
convergence of the cocycle in the sense of Definition 3.1. As a
strong limit of unitaries, $\hat{w} := \lim_{n \to \infty} \hat{u}^*_t$
is an isometry. Further we have for all unit vectors
$\hat{\xi} \in \mathcal{H} \otimes \mathcal{F}$
\[
\| q\; \hat{u}^*_t\;\hat{\xi} \|
\;\stackrel{t \to \infty}{\longrightarrow}\; 1.
\]
This has been shown for $\hat{\xi} = \xi \otimes \Omega$ before and
in general it follows by the approximations above. It implies
that the range of $\hat{w}$ is contained in $\mathcal{F}$.

Finally, we have to show the opposite inclusion, i.e., $\mathcal{F}$
is contained in the range of~$\hat{w}$. For this, note that
linear combinations of vectors of the form $\delta \otimes \eta_t$
with $\eta_t \in \mathcal{F}_{t]}$ for all~$t$, are dense in
$\delta \otimes \mathcal{F} \simeq \mathcal{F}$. From
$\hat{u}_s (\delta \otimes \Omega) = \delta \otimes \Omega$ we get
\[
\hat{u}_{t+s}(\delta \otimes \eta_t)
= \hat{u}_t (Id \otimes \gamma_t)(\hat{u}_s)\;(\delta \otimes \eta_t)
= \hat{u}_t (\delta \otimes \eta_t).
\]
We infer that for all $\eta \in \mathcal{F}$ also
$(\hat{u}_t\eta)$ is convergent. Thus our proof can be finished by
applying Lemma 3.3.
\endproof

\numeryfalse

\end{document}